\documentclass[oneside, 12pt]{amsart}

\usepackage[utf8]{inputenc}
\usepackage[margin=1in]{geometry}
\usepackage{amsthm}
\usepackage{enumitem}  
\usepackage[OT2,T1]{fontenc}
\usepackage{amscd,amssymb, enumitem, amsmath,mathrsfs, amsfonts, amsaddr}
\usepackage{url}
\usepackage{fullpage}
\usepackage{microtype}
\usepackage[backref=page]{hyperref}
\usepackage{hyperref}

\DeclareSymbolFont{cyrletters}{OT2}{wncyr}{m}{n}
\DeclareMathSymbol{\Sha}{\mathalpha}{cyrletters}{"58}

\newtheorem{theorem}{Theorem}[section]
\newtheorem{lemma}[theorem]{Lemma}
\newtheorem{proposition}[theorem]{Proposition}
\newtheorem{corollary}[theorem]{Corollary}

\theoremstyle{definition}

\newtheorem{example}[theorem]{Example}

\newtheorem{claim}[theorem]{Claim}
\newtheorem{remark}[theorem]{Remark}
\newtheorem{emptyremark}[theorem]{}

\newtheorem*{acknowledgement}{Acknowledgement}

\theoremstyle{remark}

\title{Purely additive reduction of abelian varieties with torsion}
\author{Mentzelos Melistas}
\address{Steklov Mathematical Institute of Russian Academy of Sciences, Moscow, Russia,\\ ~~~~~~~~~~~~~~~~~~~~~~~~~~~~~~~~~~ email: mentzmel@gmail.com}
\date{\today.}

\begin{document}

\maketitle

\begin{abstract}
    Let $\mathcal{O}_K$ be a discrete valuation ring with fraction field $K$ of characteristic $0$ and algebraically closed residue field $k$ of characteristic $p > 0$. Let $A/K$ be an abelian variety of dimension $g$ with a $K$-rational point of order $p$. In this article, we are interested in the reduction properties that $A/K$ can have. After discussing the general case, we specialize to $g=1$, and we study the possible Kodaira types that can occur.
\end{abstract}


\section{Introduction}

Let $\mathcal{O}_K$ be a discrete valuation ring with valuation $v_K$, fraction field $K$ of characteristic $0$ and algebraically closed residue field $k$  of characteristic $p > 0$. Let $A/K$ be an abelian variety of dimension $g$. In this paper we investigate if the existence of a $K$-rational point of order $p$ on $A/K$ imposes restrictions on the reduction properties of $A/K$. We prove Theorem \ref{thm2intro} below, using ideas from \cite{cx}.
\begin{theorem}\label{thm2intro}
Let $\mathcal{O}_K$ be a discrete valuation ring with valuation $v_K$, fraction field $K$ of characteristic $0$ and residue field $k$ which is assumed to be algebraically closed of characteristic $p > 0$. Let $A/K$ be an abelian variety with a $K$-rational point of order $p$. Let $L/K$ be an extension of minimal degree over which $A/K$ acquires semi-stable reduction. If $v_K(p)<\frac{p-1}{[L:K]} $, then $A/K$ cannot have purely additive reduction.
\end{theorem}

We then turn our attention to the case where $g=1$ and char$(k) \geq 5$. Suppose that $E/K$ is an elliptic curve and assume that $E/K$ has a $K$-rational point of order $N \geq 5$. We study the possible reduction types that $E/K$ can have. The situation where $p \nmid N$ is well-known. In this case, $E/K$ can only have semi-stable reduction (see Corollary $6.5$ of \cite{silverbergzarhin95} or Theorem $2$ of \cite{freysomeremarks}). When $p \mid N$, additive reduction for $E/K$ can also occur. The first to take up the study of elliptic curves with additive reduction and a $K$-rational point of order $p$ were Lenstra and Oort in \cite{lenstraoortpurelyadditivereduction}. We study in this article what are the possible Kodaira types of reduction that can occur in this situation. More precisely, we prove the following theorems.
\begin{theorem}\label{thm1intro}
Assume that char$(k)=p\geq5$ and let $E/K$ be an elliptic curve. Suppose also that $E/K$ has a $K$-rational point of order $p^n$ for some $n \geq 1$.
 \begin{enumerate}[topsep=2pt,label=(\roman*)]
\itemsep0em 
     \item If $v_K(p)<\frac{p^{n-1}(p-1)}{6}$, then $E/K$ can only have semi-stable reduction.
\end{enumerate}
\noindent Assume now that $E/K$ does not have semi-stable reduction.
\begin{enumerate}[topsep=2pt,,label=(\roman*),resume]
     \itemsep0em 
     \item If $v_K(p) \in [\frac{p^{n-1}(p-1)}{6}, \frac{p^{n-1}(p-1)}{4})$, then $E/K$ can only have reduction of type II.
      \item If $v_K(p)\in [\frac{p^{n-1}(p-1)}{4}, \frac{p^{n-1}(p-1)}{3})$, then $E/K$ can only have reduction of type II or III.

       \item If $v_K(p)\in [\frac{p^{n-1}(p-1)}{3}, \frac{p^{n-1}(p-1)}{2})$, then $E/K$ can only have reduction of type II, III, or IV.
\end{enumerate}
 \end{theorem}
 
 \begin{theorem}\label{theoremII*III*IV*}
  Assume that char$(k)=p\geq5$ and let $E/K$ be an elliptic curve. Suppose also that $E/K$ has a $K$-rational point of order $p^n$ for some $n \geq 1$.
  
  \begin{enumerate}[topsep=2pt, label=(\roman*)]
\itemsep0em 
\item If $v_K(p)< \frac{5p^{n-1}(p-1)}{6}$, then $E/K$ cannot have reduction of type II$^*$.
\item If $v_K(p)< \frac{3p^{n-1}(p-1)}{4}$, then $E/K$ cannot have reduction of type II$^*$ or III$^*$.
\item If $v_K(p)< \frac{2p^{n-1}(p-1)}{3}$, then $E/K$ cannot have reduction of type II$^*$, III$^*$, or IV$^*$.
\end{enumerate}
 \end{theorem}

This paper is organized as follows. Theorems \ref{thm2intro}, \ref{thm1intro}, and \ref{theoremII*III*IV*} are the main results of Section \ref{section2}. In Section \ref{section3}, we present some examples which show that the bounds of Theorems \ref{thm2intro} and \ref{thm1intro} are sharp. Section \ref{section4} is devoted to an alternative way of proving Theorems \ref{thm1intro} and \ref{theoremII*III*IV*} for small primes using explicit equations for modular curves. Finally, in the last section we explain a difference, for dimension $1$, between the equicharacteristic case, i.e., char$(K)=p$, which is studied in \cite{liedtkeschroer}, and our case. 
\begin{acknowledgement}
A large part of this work is contained in my doctoral dissertation and I gratefully acknowledge support from the University of Georgia. I would like to thank my advisor, Dino Lorenzini, for providing some of the examples that appear in Section \ref{section3} and for the many useful suggestions during the preparation of this work. I would also like to thank the anonymous referees for many insightful suggestions which improved the manuscript and for providing some arguments that were used in the proof of Theorem \ref{theoremII*III*IV*}. During the final stages of this work I was supported by the Steklov International Mathematical Center and the Ministry of Science and Higher Education of the Russian Federation (agreement no. 075-15-2019-1614).
\end{acknowledgement}

\section{Proofs of Theorems \ref{thm2intro}, \ref{thm1intro}, and $\ref{theoremII*III*IV*}$ }\label{section2}

In this section, we prove Theorems \ref{thm2intro}, \ref{thm1intro}, and \ref{theoremII*III*IV*}. Before we begin the proof we need to set up some notation as well as recall some basic facts concerning reduction of abelian varieties. Let $\mathcal{O}_K$ be a discrete valuation ring with valuation $v_K$, fraction field $K$ of characteristic $0$ and algebraically closed residue field $k$  of characteristic $p > 0$.  If $L/K$ is a finite field extension and $\mathcal{O}_K$ is complete, then we will denote by $\mathcal{O}_L$ the integral closure of $\mathcal{O}_K$ in $L$ and by $v_L$ the associated normalized discrete valuation of $L$. This notation will be fixed throughout this section.

\begin{emptyremark}\label{recallandnotation}
Let $A/K$ be an abelian variety of dimension $g$. We denote by $\mathcal{A}/\mathcal{O}_K$ the N\'eron model of $A/K$ (see \cite{neronmodelsbook} for the definition as well as the basic properties of N\'eron models). The special fiber $\mathcal{A}_k/k$ of $\mathcal{A}/\mathcal{O}_K$ is a smooth commutative group scheme. We denote by $\mathcal{A}^0_k/k$ the connected component of the identity of $\mathcal{A}_k/k$. The finite \'etale group scheme defined by $\Phi_k:=\mathcal{A}_k/\mathcal{A}^0_k$ is called the component group of $\mathcal{A}/\mathcal{O}_K$. By a theorem of Chevalley (see Theorem $1.1$ of \cite{con}) we have a short exact sequence $$0\longrightarrow T \times U \longrightarrow \mathcal{A}^0_k \longrightarrow B \longrightarrow 0, $$
where $T/k$ is a torus, $U/k$ is a unipotent group, and $B/k$ is an abelian variety. The number $\text{dim}(U)$ (resp. $\text{dim}(T)$, $\text{dim}(B)$) is called the unipotent (resp. toric, abelian) rank of $A/K$. By construction, $g=\text{dim}(U)+\text{dim}(T)+\text{dim}(B)$. We say that $A/K$ has purely additive reduction if $g=\text{dim}(U)$, or equivalently, if $\text{dim}(T)=\text{dim}(B)=0$.

Let $\mathcal{O}_K \subset \mathcal{O}_{K'}$ be a local extension of discrete valuation rings with fraction fields $K$ and $ K'$. Assume that $[K':K]$ is finite and let $e(K'/K)$ be the ramification index of $\mathcal{O}_{K'}$ over $\mathcal{O}_{K}$. We will say that the extension $K'/K$ is tame if $e(K'/K)$ is coprime to the characteristic of $k$ (note that we assume that $k$ is algebraically closed).
\end{emptyremark}

We are now ready to proceed to the proof of Theorem \ref{thm2intro}.

\begin{proof}[\it Proof of Theorem \ref{thm2intro}]

Assume that $A/K$ has purely additive reduction and that $v_K(p)<\frac{p-1}{[L:K]}$, and we will find a contradiction. Let $\mathcal{O}_{K'}$ be the completion of $\mathcal{O}_K$ and let $K'$ be the field of fractions of $\mathcal{O}_{K'}$. Since $\mathcal{O}_{K'}$ has ramification index $1$ over $\mathcal{O}_K$, the base change $A_{K'}/K'$ of $A/K$ has the same reduction as $A/K$. Therefore, we can assume that $\mathcal{O}_K$ is complete.

Recall that $L/K$ is an extension of minimal degree such that the base change $A_L/L$ has semi-stable reduction and set $m:=[L:K]$. If $L/K$ is not tame, then $p \leq m$ so $\frac{p-1}{m}<1$. Therefore, we can assume that $L/K$ is tame. We denote by $\mathcal{O}_L$ the integral closure of $\mathcal{O}_K$ in $L$, which is again a discrete valuation ring because $\mathcal{O}_K$ is complete. We denote by $v_L$ the associated (normalized) valuation of $\mathcal{O}_L$. Since $\mathcal{O}_L$ has ramification index $m$ over $\mathcal{O}_K$, the restriction $v_L|_K$ of $v_L$ to $K$ satisfies $v_L|_K=mv_K$. Let $\mathcal{A}/\mathcal{O}_K$ be the N\'eron model of $A/K$ and let $\mathcal{A}'/\mathcal{O}_L$ be the N\'eron model of $A_L/L$. If $\mathcal{A}_L/\mathcal{O}_L$ is the base change of $\mathcal{A}/\mathcal{O}_K$ to $\mathcal{O}_L$ we obtain, using the universal property of $\mathcal{A}'/\mathcal{O}_L$, an $\mathcal{O}_L$-morphism $\phi : \mathcal{A}_L \rightarrow \mathcal{A}'$ which gives rise to a $k$-morphism $\psi : (\mathcal{A}_L^0)_{k} \rightarrow (\mathcal{A'}^0)_{k}$ between the connected components of the special fibers. Since $A/K$ has purely additive reduction, 
$(\mathcal{A}_L^0)_{k}$ is a unipotent group and since $A_L/L$ has semi-stable reduction, $(\mathcal{A'}^0)_{k}$ is an extension of a torus by an abelian variety. Since there are no nonconstant morphisms from a unipotent group to an abelian variety (see Lemma $2.3$ of \cite{con}), and there are no nonconstant morphisms from a unipotent group to a torus (see Corollary 14.18 of \cite{milnealgebraicgroups}), we obtain that $\psi$ is constant.

 Recall that we assume that $A/K$ has purely additive reduction and, hence, the toric rank of $A/K$ is zero. Proposition $1.8$ of \cite{lorenziniliu} (see also \cite{complementto18} for a complement to this Proposition $1.8$) tells us that the component group of $A/K$ is killed by $[L:K]^2$ and since $L/K$ is tame, we obtain that $p$ cannot divide the order of the component group of $A/K$. Therefore, the image of $P$ under the reduction map belongs to the connected component of the identity.

  Let ${A_L}^1(L)$ denote the kernel of the reduction map $A_L(L) \rightarrow \mathcal{A'}_k(k)$. Since $v_L(p)<m(\frac{p-1}{m})=p-1 $, the reduction map of the base change is injective on torsion points (see the Appendix of \cite{katzgaloispropertiesoftorsionpts}) and, hence, $P \notin {A_L}^1(L)$. Therefore, the image of $P$ under the reduction map belongs to the connected component of the identity and has order $p$. However, this is a contradiction since $\psi$ is constant.
\end{proof}

\begin{corollary}\label{corollarydimension2}
Let $\mathcal{O}_K$ be a discrete valuation ring with valuation $v_K$, fraction field $K$ of characteristic $0$ and residue field $k$ which is assumed to be algebraically closed of characteristic $p > 0$. Let $A/K$ be an abelian surface with a $K$-rational point of order $p$. Assume that $v_K(p)=1$.
 \begin{enumerate}[topsep=2pt,label=(\roman*)]
\itemsep0em 
     \item  If $p > 13$, then $A/K$ cannot have purely additive reduction.
     \item If $p=11$ or $13$, and $A/K$ has purely additive reduction, then $A/K$ has potentially good reduction.
\end{enumerate}

\end{corollary}
\begin{proof}
We can without loss of generality assume that $\mathcal{O}_K$ is complete. Let $L/K$ be the minimal Galois extension of $K$ such that the base change $A_L/L$ of $A/K$ to $L$ has semi-stable reduction. If $x=p^{a_1}_1....p^{a_n}_n$ is a positive integer with $p_i$ distinct prime numbers, then we let $$L(x):= \begin{cases}      
      \sum_{i=1}^n p^{a_i-1}_i(p_i-1) & \text{if ord}_2(x)\neq1, \\
      L(\frac{x}{2}) & \text{if ord}_2(x)=1,  \\
      0  & \text{if } x=0,1.
   \end{cases}
  $$ 
  
\textit{Proof of $(i)$:} Assume that $A/K$ has purely additive reduction and we will arrive at a contradiction. Let $\epsilon$ be the exponent of the group $\rm{Gal}(L/K)$ and write $\epsilon=p^w \cdot \epsilon^{(p)}$, for some $w \geq 0$ with $p \nmid \epsilon^{(p)}$. If $t_L$ be the toric rank of $A_L/L$, then Proposition $3.1$ of \cite{lorenzinigroupofcomponentsofjacobians1990} tells us that $\rm{max} ( L(p^w), L(\epsilon^{(p)}) ) \leq 2\cdot2+0-t_L\leq 4$. If $w>0$, then $L(p^w) >4$ because $p>13$. Therefore, $w=0$ and, hence, $p \nmid \epsilon$. Since the exponent and the order of a finite group have the same prime divisors, we obtain that $L/K$ is tame. Since $L/K$ is also totally ramified, we have that $L/K$ is cyclic and, hence, $\epsilon=[L:K]$. Moreover, since $L(\epsilon) \leq 4$, we obtain that $\epsilon \leq 12$ by the definition of $L(\epsilon)$. Finally, since $p>13$, we find that $v_K(p)=1<\frac{p-1}{\epsilon}$ and, therefore, Theorem \ref{thm2intro} gives that $A/K$ cannot have purely additive reduction.

\textit{Proof of $(ii)$:} Assume that $p=11$ or $13$, that $A/K$ has purely additive reduction, and $A_L/L$ does not have potentially good reduction, and we will find a contradiction. Let $t_L \geq 1$ be the toric rank of $A_L/L$ and let $m=[L:K]$. Proceeding in exactly the same way as in $(i)$ we can show that $L/K$ is tame and that $L(m) \leq 2\cdot 2+0-t_L=4-t_L<4$. Therefore, $m \leq 8$ by the definition of $L(m)$. Since $p=11$ or $13$, we find that $v_K(p)=1<\frac{p-1}{m}$. Therefore, Theorem \ref{thm2intro} implies that $A/K$ cannot have purely additive reduction, which is a contradiction. This concludes the proof.
\end{proof}

\begin{remark}
Keeping the same hypotheses as in Corollary \ref{corollarydimension2}, Penniston has proved (see Theorem 6.2 of \cite{penniston2000}) that if $p>5$ and $A/K$ is the Jacobian of a smooth, projective, and geometrically connected curve such that $A(K)$ contains a point of order $p^2$, then $A/K$ cannot have purely additive reduction. Examples of abelian surfaces satisfying the hypotheses of Part $(ii)$ of Corollary \ref{corollarydimension2} can be found in Remark \ref{remarkdim2}.
\end{remark}

We now turn our attention to $g=1$. The key ingredient for the proof of Theorem \ref{thm1intro} is Proposition \ref{remarkbasechange} below. This proposition is well known to the experts, and we include a proof here for completion.

\begin{proposition}\label{remarkbasechange}
  Let $\mathcal{O}_K$ be a complete discrete valuation ring with valuation $v_K$, fraction field $K$ and algebraically closed residue field $k$ of characteristic $p \geq 5$. Let $E/K$ be an elliptic curve. Then there is an extension $L/K$ of minimal degree over which the base change $E_L/L$ of $E/K$ has
  semi-stable reduction. The degree $[L:K]$ is determined by the Kodaira type of $E/K$ as follows:
$ \begin{cases}      
      [L:K]=2 & \text{if } E/K \text{ has reduction of type I}_n^*, \text{for some }n \geq 0 .\\
      [L:K]=3 & \text{if } E/K \text{ has reduction of type IV or IV}^* .\\
      [L:K]=4 & \text{if } E/K \text{ has reduction of type III or III}^*. \\
      [L:K]=6 & \text{if } E/K \text{ has reduction of type II or II}^* .\\
   \end{cases}
  $ 
\end{proposition}
\begin{proof}
 If the reduction of $E/K$ is of type I$_n^*$ for some $n \geq 1$, then there is an extension $L/K$ over which the base change $E_L/L$ of $E/K$ has multiplicative reduction. Moreover, $[L:K]=2$ by Part $(d)$ of Theorem $14.1$ in Appendix C.14 of \cite{aec}.
 
  Assume now that $E/K$ has potentially good reduction. Let $v_K(\Delta)$ be the minimal discriminant valuation of $E/K$. Then $L$ is obtained by adjoining to $K$ the $d$-th power of a uniformizer, where $d$ is the minimal positive integer such that $dv_K(\Delta)$ is divisible by $12$, as we now explain. Let $\mathcal{O}_L$ be the integral closure of $\mathcal{O}_K$ in $L$. We also denote by $v_L$ the associated (normalized) valuation of $\mathcal{O}_L$. Then $v_L(\Delta)=dv_K(\Delta)$. We know that the valuation of minimal discriminant for the base $E_L/L$ must be between $0$ and $12$ by Remark VII.1.1 of \cite{aec}. Moreover, when we perform a change of Weierstrass equation for an elliptic curve the valuation of the discriminant changes by a multiple of $12$. Therefore, since $12 \mid v_L(\Delta)$, $E_L/L$ has good reduction. The corresponding values for $d$ depending on each Kodaira type, which can be found using the table on page 46 of \cite{tatealgorithm}, give the degree of $L/K$. This concludes the proof of our proposition.
\end{proof}
\begin{proof}[\it Proof of Theorem \ref{thm1intro}]
By an argument similar to the start of the proof of Theorem \ref{thm2intro} we can assume that $\mathcal{O}_K$ is complete.

Let $\widetilde{E}_{ns}(k)$ denote the set of non-singular $k$-rational points of $\widetilde{E}/k$, $E^1(K)$ denote the kernel of the reduction map, and $E^0(K)$ denote the set of points whose image under the reduction map is a non-singular.
 
We first prove that if $v_K(p)<\frac{p^{n-1}(p-1)}{2}$, then it is not possible for $E/K$ to have reduction of type I$_m^*$ for any $m \geq 0$. Assume that $E/K$ has reduction of type I$_m^*$ for some $m \geq 0$ and we will find a contradiction. Let $P \in E(K)[p^n]$ be the $K$-rational point of order $p^n$. Since $E/K$ has additive reduction, we obtain that $|E(K)/E^0(K)| \leq 4$. Therefore, since $p \geq 5$ and $P$ has order $p^n$, we get that $P \in E^0(K)$.

 Assume first that $m> 0$. Let $L/K$ be an extension of $K$ over which $E_L/L$ achieves semi-stable reduction. Assume also that the degree $[L:K]$ is minimal. By Proposition \ref{remarkbasechange}, we obtain that $[L:K]=2$. Let $\mathcal{O}_L$ be the integral closure of $\mathcal{O}_K$ in $L$, which is again a discrete valuation ring because $\mathcal{O}_K$ is complete. We also denote by $v_L$ the associated (normalized) valuation of $\mathcal{O}_L$. Since $\mathcal{O}_L$ has ramification index $2$ over $\mathcal{O}_K$, we have that the restriction $v_L|_K$ of $v_L$ to $K$ satisfies $v_L|_K=2v_K$. Using Theorem $5.9.4$ of \cite{Krummthesis}, since $v_L(p)=2v_K(p)< 2(\frac{p^{n-1}(p-1)}{2})=p^{n-1}(p-1)$, we obtain that $P \not \in {E_L}^1(L)$. Let $\mathcal{E}/\mathcal{O}_K$ be the N\'eron model of $E/K$ and let $\mathcal{E}'/\mathcal{O}_L$ be N\'eron model of $E_L/L$. By looking at the base change $\mathcal{E}_L/\mathcal{O}_L$ of $\mathcal{E}/\mathcal{O}_K$ over $\mathcal{O}_L$ we obtain, using the universal property of $\mathcal{E}'/\mathcal{O}_L$, that there is an $\mathcal{O}_L$-morphism $\phi : \mathcal{E}_L \rightarrow \mathcal{E}'$ which gives rise to a $k$-morphism $(\mathcal{E}_L^0)_{k} \rightarrow (\mathcal{E}'^0)_{k}$. Since $E/K$ has additive reduction, we get that $(\mathcal{E}_L^0)_{k} \cong \mathbb{G}_{a,k}$ and on the other hand since $E_L/L$ has multiplicative reduction, we obtain that $(\mathcal{E}'^0)_{k} \cong \mathbb{G}_{m,k}$. However, because of the fact that there are no non-constant $k$-morphisms from $\mathbb{G}_{a,k}$ to $\mathbb{G}_{m,k}$  (see Corollary 14.18 of \cite{milnealgebraicgroups}) we get that the image  of $(\mathcal{E}_L^0)_{k}$ is the identity of $(\mathcal{E}'^0)_{k}$. This is a contradiction since this means that $P \in {E_L}^1(L)$ and this contradicts the fact that $P \not \in {E_L}^1(L)$ showed above. This proves that if $v_K(p)<\frac{p^{n-1}(p-1)}{2}$, then it is not possible for $E/K$ to have reduction of type I$_m^*$ for any $m > 0$. The case where $m=0$ follows directly from Theorem \ref{thm2intro}, using Proposition \ref{remarkbasechange}.
 
Using similar arguments as the proof above combined with Proposition \ref{remarkbasechange} we can prove the following three statements. If $v_K(p)<\frac{p^{n-1}(p-1)}{3}$, then it is not possible for $E/K$ to have reduction of type IV or IV$^*$. If $v_K(p)<\frac{p^{n-1}(p-1)}{4}$, then it is not possible for $E/K$ to have reduction of type III or III$^*$. If $v_K(p)<\frac{p^{n-1}(p-1)}{6}$, then it is not possible for $E/K$ to have reduction of type II or II$^*$. 

Finally, using Theorem \ref{theoremII*III*IV*} we can rule out the cases II$^*$, III$^*$, and IV$^*$ in each part of the theorem. Putting now everything together, we have proved Theorem \ref{thm1intro}.
\end{proof}

\begin{proof}[\it Proof of Theorem \ref{theoremII*III*IV*}]
By an argument similar to the start of the proof of Theorem \ref{thm2intro} we can assume that $\mathcal{O}_K$ is complete. Let $E/K$ be an elliptic curve with additive and potentially good reduction. Let $L/K$ be an extension of $K$ of minimal degree over which the base change $E_L/L$ of $E/K$ to $L$ acquires good reduction. Let $\mathcal{O}_L$ be the integral closure of $\mathcal{O}_K$ in $L$, which is again a discrete valuation ring because $\mathcal{O}_K$ is complete. We also denote by $v_L$ the associated (normalized) valuation of $\mathcal{O}_L$. Note that the restriction $v_L|_K$ of $v_L$ to $K$ satisfies $v_L|_K=mv_K$, where $m:=[L:K]$.

Since we assume that $p \geq 5$, we can find short  minimal Weierstrass equations
\begin{align*}
    W_{min}: \; y^2=x^3+ax+b \; \text{  and   } \; W_{min}^L: \; y^2=x^3+a'x+b'
\end{align*}
for $E/K$ and $E_L/L$ respectively. We will denote by $\Delta_{min}$ the discriminant of $W_{min}$.

\begin{lemma}\label{lemmavaluation}
Let $P \in E(K)$ such that $P \not\in E^1(K)$ but $P \in {E_L}^1(L)$. If $P$ corresponds to $(x_0,y_0)$ on $W_{min}$, then $v_L(x_0)=0=v_L(y_0)$. Moreover, the corresponding point $(x_0',y_0')$ on $W_{min}^L$ satisfies 
\begin{align*}
    v_L(x_0')=\frac{-v_L(\Delta_{min})}{6} \; \text{  and   } \; v_L(y_0')=\frac{-v_L(\Delta_{min})}{4}.
\end{align*}
\end{lemma}
\begin{proof}
 There exists an isomorphism $W_{min}\longrightarrow W_{min}^L$ defined over $L$ which is of the form $(x,y) \mapsto (u^{-2}x,u^{-3}y)$, for some $u \in L^{*}$. Since this isomorphism alters the discriminant by $u^{-12}$ and $E_L/L$ has good reduction by assumption, we find that $v_L(u)=\frac{v_L(\Delta_{min})}{12}.$
 
 Since $(x_0', y_0'):=(u^{-2}x_0, u^{-3}y_0) \in {E_L}^1(L)$ and, hence, $(x_0', y_0')$ reduces to the point at infinity we have 
 \begin{align*}
     v_L(u^{-2}x_0)<0, \; v_L(u^{-3}y_0)<0 \text{ and } 3v_L(u^{-2}x_0)=2v_L(u^{-3}y_0).
 \end{align*}
 Therefore, $v_L(x_0)< 2v_L(u)=\frac{v_L(\Delta_{min})}{6}$ and $v_L(x_0)<3v_L(u)=\frac{v_L(\Delta_{min})}{6}$.
 
 On the other hand, since $P \not\in E^1(K)$ we obtain that $v_K(x_0)\geq 0$ and $v_K(y_0) \geq 0$. Moreover, $3v_L(u^{-2}x_0)=2v_L(u^{-3}y_0)$ implies that $3v_L(x_0)=2v_L(y_0)$ and, hence, $3v_K(x_0)=2v_K(y_0)$. Write $v_K(x_0)=2r$, for some $r \geq 0$. Since $E/K$ has potentially good reduction, we have that $v_K(\Delta_{min})<12$ (see Table 4.1 of \cite{silverman2}). This implies that $v_K(x_0) < 2$ and since $v_K(x_0)$ is even we obtain that $v_K(x_0)=0$. Therefore, $v_L(x_0)=0$ and also $v_L(y_0)=0$ because $3v_L(x_0)=2v_L(y_0)$.
 
 Finally, since $v_L(u)=\frac{v_L(\Delta_{min})}{12}$ and $(x_0', y_0')=(u^{-2}x_0, u^{-3}y_0)$, we obtain that $v_L(x_0')=\frac{-v_L(\Delta_{min})}{6}$ and  $v_L(y_0')=\frac{-v_L(\Delta_{min})}{4}$, as needed. This proves the lemma.

\end{proof}

\begin{corollary}\label{corollarypointorder}
Suppose that $v_L(p)< \frac{(p^n-p^{n-1})v_L(\Delta_{min})}{12}$ for some $n \geq 1$ and that $P\in E(K)$ has order $p^n$. Then $P \not\in {E_L}^1(L)$.
\end{corollary}

\begin{proof}
 Suppose that $P \in {E_L}^1(L)$ and we will find a contradiction. If $P$ corresponds to $(x_0',y_0')$ on $W_{min}^L$, then Lemma \ref{lemmavaluation} gives that $v_L(\frac{-x_0'}{y_0'})=\frac{v_L(\Delta_{min})}{12}$. On the other hand, Theorem IV.6.1 of \cite{aec} gives that $v_L(\frac{-x_0'}{y_0'}) \leq \frac{v_L(p)}{p^n-p^{n-1}}< \frac{v_L(\Delta_{min})}{12}$. This is a contradiction and we have proved our corollary.
\end{proof}

We are now ready to conclude our proof.

\textit{Proof of $(i)$:} Assume that $v_K(p)< \frac{5p^{n-1}(p-1)}{6}$ and that $E/K$ has reduction of type II$^*$ and we will arrive at a contradiction. Let $P$ a point of order $p^n$. Proceeding in a similar way as in the proof of Theorem \ref{thm2intro} we will find a contradiction provided that we can show that $P \not\in {E_L}^1(L)$. Since $E/K$ has reduction of type II$^*$, Proposition \ref{remarkbasechange} gives that $[L:K]=6$. Moreover, from Table 4.1 of \cite{aec} we obtain that $v_K(\Delta_{min})=10$. Therefore, $v_L(\Delta_{min})=60$ and $v_L(p)=6v_K(p)$. This implies that if $v_K(p)< \frac{5p^{n-1}(p-1)}{6}$, then $v_L(p)< \frac{(p^n-p^{n-1})v_L(\Delta_{min})}{12}$ and, hence, Corollary \ref{corollarypointorder} gives that $P \not\in {E_L}^1(L)$, as needed. This proves $(i)$.

The proofs of parts $(ii)$ and $(iii)$ are similar using Proposition \ref{remarkbasechange} and Table 4.1 of \cite{silverman2}.

\end{proof}

Our next proposition shows that the statements of Theorem \ref{theoremII*III*IV*} are the best that one can hope for.

\begin{proposition}\label{prop:basechange}
Let $\mathcal{O}_K$ be a complete discrete valuation ring with valuation $v_K$, fraction field $K$ of characteristic $0$ and residue field $k$ which is assumed to be algebraically closed of characteristic $p \geq 5$. Let $E/K$ be an elliptic curve with a rational point of order $p \geq 5$ and assume that $E/K$ has additive reduction.
 \begin{enumerate}[topsep=2pt,label=(\roman*)]
\itemsep0em 
  \item If $E/K$ has reduction of type II and $L/K$ is a finite extension with $v_L(p)=5v_K(p)$, then $E_L/L$ has reduction of type II$^*$.
 \item If $E/K$ has reduction of type III and $L/K$ is a finite extension with $v_L(p)=3v_K(p)$, then $E_L/L$ has reduction of type III$^*$.
 \item If $E/K$ has reduction of type IV and $L/K$ is a finite extension with $v_L(p)=2v_K(p)$, then $E_L/L$ has reduction of type IV$^*$.
 \end{enumerate}
 \end{proposition}
\begin{proof}
The proof is a consequence of the table on page $46$ of \cite{tatealgorithm}. Assume $E/K$ has reduction of type II (resp. III, IV) and let $\Delta$ be the minimal discriminant of $E/K$. By the table on page $46$ of \cite{tatealgorithm} $v_K(\Delta)=2$ (resp. $v_K(\Delta)=3$,$v_K(\Delta)=4$). Let $L/K$ be a field extension with $v_L(p)=5v_K(p)$ (resp. $v_L(p)=3v_K(p)$, $v_L(p)=2v_K(p)$). Then $v_L(\Delta)=10$ (resp. $v_L(\Delta)=9$, $v_L(\Delta)=8$) so by the table on page $46$ of \cite{tatealgorithm} we get that the base change of $E/K$ to $L$, $E_L/L$, has reduction type of type II$^*$ (resp. III$^*$, IV$^*$). This proves the proposition.
\end{proof}

\begin{remark}
When $p=3$ and $v_K(p)=1$, Kozuma (see Proposition 3.5 and Lemma 3.6 of \cite{koz}) has described the possible reduction types of elliptic curves $E/K$ that have a $K$-rational point of order $3$ using explicit Weierstrass equations. We note that in this case the reduction types II$^*$ and III$^*$ cannot occur.
\end{remark}

\section{Some Examples}\label{section3}

In this section, we present examples of abelian varieties with torsion points and purely additive reduction. Our examples show that the ramification bounds in Theorems \ref{thm2intro} and \ref{thm1intro} are sharp.

\begin{example}
Let $p$ be an odd prime and $s$ be an integer with $1 \leq s \leq p-2$. Consider the smooth projective curve $C_{p,s}/\mathbb{Q}$ birational to $$y^p=x^s(1-x).$$
The curve $C_{p,s}/\mathbb{Q}$ has genus $\frac{p-1}{2}$. The Jacobian $J_{p,s}/\mathbb{Q}$ of $C_{p,s}/\mathbb{Q}$ has a $\mathbb{Q}$-rational point of order $p$ (see Theorem $1.1$ of \cite{someresultsgrossrohrlich}).
Each $C_{p,s}/\mathbb{Q}$ is a quotient of the Fermat curve $F_p/\mathbb{Q}$ (see Example 5.1 of \cite{lorenzini1993} for more details). It turns out that each $J_{p,s}/\mathbb{Q}$ has good reduction away from $p$ and potentially good reduction modulo $p$. When  $C_{p,s}/\mathbb{Q}$ is tame (see Example 5.1 of \cite{lorenzini1993} for the definition) $J_{p,s}/\mathbb{Q}$ has purely additive reduction modulo $p$ and achieves good reduction after a totally ramified extension of degree $2(p-1)$ (see Anmerkung in page 339 of \cite{maedastablereductionfermatcurve} for the last statement).
\end{example}

%\begin{example}
%Let $J_0(13^3)/\mathbb{Q}$ be the Jacobian of the modular curve $X_0(13^3)/\mathbb{Q}$. The variety $J_0(13^3)/\mathbb{Q}$ has a $\mathbb{Q}$-rational point of order $13$ (see Theorem 1 of \cite{lingqrationalcuspidalgroup}) and has good reduction away from $13$. Moreover, $J_0(13^3)/\mathbb{Q}$ has potentially good reduction modulo $13$ (see Theorem 9.4 of \cite{colemanmcmurdy}) and has purely additive reduction modulo $13$ (see Section 1.4 of \cite{edixhovenminimalresolutionofx0(n)} and note that the genus of $X_0(13)/\mathbb{Q}$ is $0$).
%\end{example}

Let $K$ be as in the previous section. Let $C/K$ be a smooth, proper, and geometrically connected curve of genus $2$ with Jacobian $J/K$.  Assume that $J(K)$ contains a point of prime order $p$. When $p=7,11,13$ and $v_K(p)=1$, purely additive reduction for $J/K$ is still allowed by Corollary \ref{corollarydimension2}. The examples below when $p=7$ show that the special fiber of the minimal regular model of the curve $C/K$ is surprisingly simple. 

\begin{example}\label{example7torsion1genus2}
Consider the curve $C/\mathbb{Q}$ given by the following equation 
$$y^2+(x^3+1)y=22x^6-10x^5-24x^4+5x^3+11x^2-x-2.$$
This is a genus $2$ curve with LMFD \cite{lmfdb} label 40817.a.40817.1. The Jacobian $\text{Jac}(C)/\mathbb{Q}$ of $C/\mathbb{Q}$ has a $\mathbb{Q}$-rational point of order $7$. 

Using the command \texttt{genus2reduction()} in SAGE \cite{sagemath} we see that $C/\mathbb{Q}$ has reduction type [VIII-1] (\cite{namikawauenoclassification}, page 156) modulo $7$ and, hence, potentially good reduction modulo $7$. The special fiber of the minimal regular model is irreducible with multiplicity $1$. Let $C_{\mathbb{Q}_7^{unr}}/\mathbb{Q}_7^{unr}$ be the base change of $C/\mathbb{Q}$ to the maximal unramified extension of $\mathbb{Q}_7$. Since the minimal strict normal crossings model of $C_{\mathbb{Q}_7^{unr}}/\mathbb{Q}_7^{unr}$ contains no irreducible components of genus greater than $0$, by Theorem $6.1$ of \cite{neronmodelslorenzini} the abelian rank of $\text{Jac}(C)/\mathbb{Q}$ modulo $7$ is zero. Therefore, since $C/\mathbb{Q}$ has potentially good reduction modulo $7$, $\text{Jac}(C)/\mathbb{Q}$ has purely additive reduction modulo $7$.
\end{example}

 One way of trying to find a curve of genus $2$ with a Jacobian having a $K$-rational point of order $p$ is to find a degree $2$ cover of an elliptic curve having a $K$-rational point of order $p$. This is our method for producing Examples \ref{exampledegree2coverjacobianV} and \ref{exampledegree2coverjacobianIIII0} below.

\begin{example}\label{exampledegree2coverjacobianV}
Consider the curve $C/\mathbb{Q}$ defined by the equation
$$y^2=x^6 - 182763x^2 + 31201254.$$
This is a curve of genus $2$. Using the command \texttt{TorsionSubgroup()} in MAGMA \cite{magmareference}, we find that the Jacobian $\text{Jac}(C)/\mathbb{Q}$ of $C/\mathbb{Q}$ has a $\mathbb{Q}$-rational point of order $7$.  

Using the command \texttt{genus2reduction()} in SAGE \cite{sagemath} we see that $C/\mathbb{Q}$ has reduction type [V] (\cite{namikawauenoclassification}, page 156) modulo $7$, and, hence, potentially good reduction. The special fiber of the minimal regular model consists of two irreducible components of multiplicity $1$. Using a similar argument as in Example \ref{example7torsion1genus2}, we find that $\text{Jac}(C)/\mathbb{Q}$ has purely additive reduction modulo $7$. Note that, by Table 1 of \cite{Liugenus2algorithm}, the minimal degree of an extension over which the base change ${\text{Jac}(C)}_{\mathbb{Q}_7^{unr}}/\mathbb{Q}_7^{unr}$ acquires good reduction is $6$. This shows that Theorem \ref{thm2intro} is sharp in the sense that in our example $\frac{p-1}{m}=1$ and $v_K(p)=1$, but ${\text{Jac}(C)}_{\mathbb{Q}_7^{unr}}/\mathbb{Q}_7^{unr}$ has purely additive reduction.
\end{example}

\begin{example}\label{exampledegree2coverjacobianIIII0}
Consider the curve $C/\mathbb{Q}$ defined by the equation
$$y^2=x^6 + \frac{1}{4}x^4 - 141x^2 + 657.$$

This is a curve of genus $2$. Using the command \texttt{TorsionSubgroup()} in MAGMA \cite{magmareference}, we find that the Jacobian $\text{Jac}(C)/\mathbb{Q}$ of $C/\mathbb{Q}$ has a $\mathbb{Q}$-rational point of order $7$. 

Using the command \texttt{genus2reduction()} in SAGE \cite{sagemath} we see that $C/\mathbb{Q}$ has reduction type [II-II-0] (\cite{namikawauenoclassification}, page 163) modulo $7$, and, hence, potentially good reduction. The special fiber of the minimal regular model is irreducible of multiplicity $1$. Therefore, $\text{Jac}(C)/\mathbb{Q}$ has potentially good reduction modulo $7$.
\end{example}

\begin{remark}
When $g=2$, $p=11$, and $v_K(p)=1$, purely additive reduction for Jacobians of dimension $g$ with a $K$-rational point of order $p$ is still allowed. In order for such an example to exist, the minimal degree of an extension over which $J/K$ acquires good reduction must be $10$ or $12$. Families of curves of genus $2$ defined over $\mathbb{Q}$ whose Jacobian has a $\mathbb{Q}$-rational point of order 11 have been constructed in \cite{flynntorsionjacobians}, in \cite{daowsudschmidt}, and some sporadic examples can be found in \cite{leprevost11torsion}. We performed a search through those families. For all the examples that we found there, the reduction of the Jacobian modulo $p=11$ had either a positive toric rank or good reduction.
\end{remark}

\begin{example}
There exists a field $K$ and a prime number $p$ as in Theorem \ref{thm1intro}, Part $(iv)$, such that $v_K(p)=\frac{p-1}{2}$ and such that there exists an elliptic curve $E_K/K$ with a $K$-rational point of order $p$ and Kodaira type I$_1^*$. Indeed, consider the curve $E/\mathbb{Q}(\sqrt{5})$ given by $$ y^2+\Big(\frac{1-\sqrt{5}}{2}+1\Big)xy+\Big( \frac{1-\sqrt{5}}{2}+1\Big)y=x^3+x^2+\Big(\frac{-314(1-\sqrt{5})}{2}-1031 \Big)x+\frac{5958(1-\sqrt{5})}{2}+12717.$$ This is the curve with label 2.2.5.1-1100.1-i1  in the LMFDB database \cite{lmfdb}. Let $M=\mathbb{Q}(\sqrt{5})$ and $\mathfrak{p}$ be the prime above $p:=5$ in $M$. If $K:=M_{\mathfrak{p}}^{unr}$ is the maximal unramified extension of the completion of $K$ with respect to $\mathfrak{p}$, then the base change $E_K/K$ of $E/\mathbb{Q}(\sqrt{5})$ has a $K$-rational torsion point of order $10$, because $E/\mathbb{Q}(\sqrt{5})$ does. Moreover, $E_K/K$ has additive reduction of type I$_1^*$. The proof of Theorem \ref{thm1intro} does not apply to this example because, using the same notation as in the proof, looking at the elliptic curve over the minimal extension of degree $2$ over $\mathbb{Q}(\sqrt{5})$ in which $E$ achieves semi-stable reduction, say $L$, we get that $2P \in (E_L)_1(L)$.
\end{example}

\begin{example}\label{exampleappendix}
There exists a field $K$ and a prime number $p$ as in Theorem \ref{thm1intro}, Part $(iv)$, such that $v_K(p)=\frac{p-1}{2}$ and such that there exists an elliptic curve $E_K/K$ with a $K$-rational point of order $p$ and Kodaira type I$_0^*$. In the tables of Van Hoeij \cite{vanhoeijdata} there is an elliptic curve $E/M$ (with $N = 58$, deg$v = 28$, deg$j = 1$, $j = -3375$ in Van Hoeij's notation) which is defined over a number field $M$ with $[M:\mathbb{Q}]=28$ and with a $M$-rational point of order $29$ (in fact it has a point of order $58$).  This elliptic curve is a CM elliptic curve and there is a prime ideal $\mathfrak{p}$ of $M$ above $p:=29$ such that $e(\mathfrak{p}|29)=\frac{29-1}{2}=14$. If $K:=M_{\mathfrak{p}}^{unr}$ is the maximal unramified extension of the completion of $K$ with respect to $\mathfrak{p}$, then the base change $E_K/K$ of $E/M$ has a $K$-rational torsion point of order $29$, because $E/M$ does. Moreover, $E_K/K$ has additive reduction of type I$_0^*$.
\end{example}

\begin{example}\label{exampleIV}
There exists a field $K$ and a prime number $p$ as in Theorem \ref{thm1intro}, Part $(iii)$, such that $v_K(p)=\frac{p-1}{3}$ and such that there exists an elliptic curve $E_K/K$ with a $K$-rational point of order $p$ and Kodaira type IV. Let $$M=\mathbb{Q}[e]/(e^{12}-4e^{11}+9e^{10}-6e^9-2e^8-3e^7+11e^6-12e^4+7e^3+2e^2-3e+1)$$ and consider the elliptic curve $E/M$ given by $$E(b,c): y^2+(1-c)xy-by=x^3-bx^2$$ where $b,c$ are given as follows. Let $g=\frac{u}{1042}$ where $$u=-(530e^{11}-1314e^{10}+2235e^9+1428e^8-646e^7-5201e^6-658e^5+4866e^4+808e^3-574e^2-1834e-99)$$  and let  $$r=\frac{e^2g-eg+g-1}{e^2g-e}, \quad s=\frac{eg-g+1}{eg}.$$ Finally, put $$b=rs(r-1) \quad \mathrm{and} \quad c=s(r-1).$$ We have that $E/M$ is an elliptic curve with an $M$-rational point of order $p:=19$ (in fact it has a point of order $57$) and $[M:\mathbb{Q}]=12$. This is a CM elliptic curve that appears in the tables of Van Hoeij \cite{vanhoeijdata}. Let $\mathfrak{p}$ be the prime ideal $(19,e+1)$ of $M$. We have that $e(\mathfrak{p}|19)=\frac{19-1}{3}=6$. If $K:=M_{\mathfrak{p}}^{unr}$ is the maximal unramified extension of the completion of $K$ with respect to $\mathfrak{p}$, then the base change $E_K/K$ of $E/M$ has a $K$-rational torsion point of order $19$, because $E/M$ does. Moreover, $E_K/K$ has additive reduction of type IV.
\end{example}

\begin{example}\label{exampleIII}
There exists a field $K$ and a prime number $p$ as in Theorem \ref{thm1intro}, Part $(ii)$, such that $v_K(p)=\frac{p-1}{4}$ and such that there exists an elliptic curve $E_K/K$ with a $K$-rational point of order $p$ and Kodaira type III. Let $$M=\mathbb{Q}[e]/(e^8+4e^7+7e^6+8e^5+8e^4+6e^3+4e^2+2e+1)$$ and consider the elliptic curve $E/M$ given by $$E(b,c): y^2+(1-c)xy-by=x^3-bx^2$$ where $$b=\frac{-11e^7-34e^6-38e^5-16e^4-11e^3-11e^2-e+3}{17}$$ and $$c=\frac{9e^7+31e^6+59e^5+60e^4+50e^3+30e^2+25e+6}{17}.$$ This elliptic curve appears  in the table of page $18$ of \cite{ccrs} and is a CM elliptic curve. We have that $[M:\mathbb{Q}]=8$ and that $E/M$ has a $M$-rational point of order $p:=17$ (in fact it even has a $M$-rational point of order $34$). Let $\mathfrak{p}$ be the prime ideal $(17,e-7)$ of $M$. Then using SAGE we can check that $e(\mathfrak{p}|17)=\frac{17-1}{4}=4$. If $K:=M_{\mathfrak{p}}^{unr}$ is the maximal unramified extension of the completion of $K$ with respect to $\mathfrak{p}$, then the base change $E_K/K$ of $E/M$ has additive reduction of type III.
\end{example}

%\begin{example}\label{exampleII}
%There exists a field $K$ and a prime number $p$ as in Theorem \ref{thm1intro}, Part $(i)$, such that $v_K(p)=\frac{p-1}{6}$ and such that there exists an elliptic curve $E_K/K$ with a $K$-rational point of order $p$ and Kodaira type II. Let $$M=\mathbb{Q}[b]/(b^4 + 4b^3+ 78b^2 + 13b + 1)$$ and consider the elliptic curve $E/M$ given by $$E(b,c): y^2+(1-c)xy-by=x^3-bx^2,$$ where $$c=\frac{16b^3+44b^2+1354b+45}{483}.$$ This example appears in the table of page $16$ of \cite{ccrs}. The elliptic curve $E/M$ has an $M$-rational point of order $p:=13$. We have that $[M:\mathbb{Q}]=4$. Let $\mathfrak{p}$ be the prime ideal $(13,b^3 37/483+b^2 142/483+b2789/483+899/483)$ of $M$. We have that $\mathfrak{p}$ is a prime above $13$ and $e(\mathfrak{p}|13)=\frac{13-1}{6}=2$. If $K:=M_{\mathfrak{p}}^{unr}$ is the maximal unramified extension of the completion of $K$ with respect to $\mathfrak{p}$, then the base change $E_K/K$ of $E/M$ has additive reduction of type II. All of these claims can be checked using SAGE.
%\end{example}

\begin{remark}\label{remarkII}
 There exists a field $K$ and a prime number $p$ as in Theorem \ref{thm1intro}, Part $(i)$, such that $v_K(p)=\frac{p-1}{6}$ and such that there exists an elliptic curve $E_K/K$ with a $K$-rational point of order $p$ and Kodaira type II. In fact, using SAGE \cite{sagemath} we found many examples of elliptic curves that are defined over number fields, have a point of order $p:=13$, and have reduction type II modulo some prime $\mathfrak{p}$ of the ring of integers $\mathcal{O}_K$ that lies above $p$ and such that the ramification index of $\mathfrak{p}$ over $p$ is equal to $\frac{p-1}{6}$. We now explain the method we used. If $K$ is any number field and $E/K$ is an elliptic curve with a $K$-rational point of order $13$, then there exist $s,t \in K$ which are solutions to the modular curve $X_1(13):s^2=t^6-2t^5+t^4-2t^3+6t^2-4t+1,$ and such that $E/K$ can be given by the following equation
 \begin{align*}
    y^2+axy+by=x^3+bx^2,
\end{align*}
where $$a=\frac{(t-1)^2(t^2+t-1)s-t^7+2t^6+3t^5-2t^4-5t^3+9t^2-5t+1}{2t^5},$$
$$b=\frac{(t-1)^2((t^5+2t^4-5t^2+4t-1)s-t^8-t^7+4t^6+2t^5+t^4-13t^3+14t^2-6t+1)}{2t^9}$$ (see also \cite{rab} page 14). Consider $t_n:=13n-3$, for $n \in \mathbb{N}$ and let $(s_n,t_n)$ be a point of $X_1(13)(\mathbb{Q}(s_n))$. Curves corresponding to the points $(s_n,t_n)$, for $n = 1,2,...,3000$, have reduction type II modulo the prime above $13$. It seems likely that curves corresponding to the points $(s_n,t_n)$ have reduction type II modulo the prime above $13$ for all $n \in \mathbb{N}$.
\end{remark}

\begin{remark}\label{remarkII11}
There exist elliptic curves $E/K$ defined over quadratic number fields $K/\mathbb{Q}$ ramified at $p:=11$ with a $K$-rational point of order $p$ and with reduction type II modulo the prime above $p$. We now explain how to find such examples. If $K$ is any number field and $E/K$ is an elliptic curve with a $K$-rational point of order $11$, then there exist $s,t \in K$ which are solutions to $X_1(11):s^2-s=t^3-t^2,$ and such that $E/K$ can be given by the following equation \begin{align*}
   y^2+(st+t-s^2)xy+s(s-1)(s-t)t^2y=x^3+s(s-1)(s-t)tx^2
\end{align*} (see \cite{rab} page 10). Consider $t_n:=11n-3$, for $n \in \mathbb{N}$ and let $(s_n,t_n)$ be a point of $X_1(11)(\mathbb{Q}(s_n))$. Curves corresponding to the points $(s_n,t_n)$, for $n = 1,2,...,3000$, have reduction type II modulo the prime above $11$. It seems likely that curves corresponding to the points $(s_n,t_n)$ have reduction type II modulo the prime above $11$ for all $n \in \mathbb{N}$.
\end{remark}

\begin{remark}\label{remark7torsion}
When $p:=7$ and $v_K(7)=1$, then there exist elliptic curves with a $K$-rational point of order $7$ and reduction type II. For instance, the elliptic curve with Cremona label $294b2$ has a $\mathbb{Q}$-rational point of order $7$ and reduction type II modulo $7$. 

When $p:=7$ and $v_K(7)=2$, then there exist elliptic curves with a $K$-rational point of order $7$ and reduction types III and IV. Indeed, if an elliptic curve $E/\mathbb{Q}$ has reduction type II modulo $7$, then after a quadratic extension $K/\mathbb{Q}$ ramified at $7$ the base change $E_K/K$ will have reduction type IV modulo the prime of $K$ above $7$. Moreover, using the LMFDB database one can find examples of elliptic curves defined over a quadratic extension $K/\mathbb{Q}$ that have $K$-rational point of order $7$ and reduction type III modulo the prime of $K$ above $7$. One such example is the elliptic curve $E/K$ with LMFDB label $2.2.21.1$-$980.1$-p$1$ defined over $\mathbb{Q}(\sqrt{21})$.
\end{remark}

\begin{remark}\label{remarkdim2}
Let $p:=11$ or $13$ and let $K/\mathbb{Q}$ be a quadratic number field ramified at $p$. Let $E/K$ be an elliptic curve with a $K$-rational point of order $p$. Assume that $E/K$ has reduction type II modulo the prime $\mathfrak{p}$ above $p$ (see Remarks \ref{remarkII} and \ref{remarkII11} for such examples). Then the Weil restriction $\text{Res}_{K/\mathbb{Q}}(E)/\mathbb{Q}$ of $E/K$ from $K$ to $\mathbb{Q}$ has a $\mathbb{Q}$-rational point of order $p$ and purely additive reduction modulo $\mathfrak{p}$, as we now explain. The base change $\text{Res}_{K/\mathbb{Q}}(E)_K/K$ of $\text{Res}_{K/\mathbb{Q}}(E)/\mathbb{Q}$ to $K$ is isomorphic over $K$ to $E\times E^{\sigma}$, where $E^{\sigma}/K$ is the base change of $E/K$ along the map $\sigma : K \rightarrow K$. Moreover, $E/K$ has additive reduction modulo $\mathfrak{p}$, and since $\mathfrak{p}$ is the unique prime above $p$, it follows that $E^{\sigma}/K$ has additive reduction modulo $\mathfrak{p}$. Therefore, the toric and abelian ranks of $\text{Res}_{K/\mathbb{Q}}(E)_K/K$ are equal to zero. Since the toric and abelian ranks cannot decrease upon base change, it follows that $\text{Res}_{K/\mathbb{Q}}(E)/\mathbb{Q}$ has purely additive reduction modulo $p$.
\end{remark}

\section{An alternative way to prove Theorems \ref{thm1intro} and \ref{theoremII*III*IV*} for small primes}\label{section4}
Let $\mathcal{O}_K$ be a discrete valuation ring with valuation $v_K$, fraction field $K$ of characteristic $0$ and algebraically closed residue field $k$ of characteristic $p \geq 5$. We now present an alternative method to prove some partial results towards Theorems \ref{thm1intro} and \ref{theoremII*III*IV*} using explicit equations for modular curves. The idea is that using explicit equations for modular curves we can write down the general equation of an elliptic curve with a torsion point of order $p$ and then we can use Tate's algorithm \cite{tatealgorithm} to analyze the reduction of the curve based on the equation that we have. We illustrate the method for $p=5$ below and we refer to the author's Ph.D. thesis \cite{mentzelosthesis} for similar results for $p \leq 23$.

\begin{proposition}\label{5torsionprop}
Assume char$(k)=5$ and let $E/K$ be an elliptic curve. Suppose also that $E/K$ has a $K$-rational torsion point of order $5$ and that $E/K$ has additive reduction. Then
 \begin{enumerate}[topsep=2pt,label=(\roman*)]
\itemsep0em 
    \item  If $v_K(5)=1$, then $E/K$ can only have reduction type II or III.
    \item If $v_K(5)=2$, then $E/K$ can only have reduction type II, III, IV or I$_n^*$ for $n \geq 0$.
    \item If $v_K(5)=3$, then $E/K$ cannot have reduction type II$^*$.
    \item If $v_K(5)=1$ and $E/K$ has a $K$-rational torsion point of order $10$, then $E/K$ can only have reduction type III.
\end{enumerate}
 \end{proposition}
 \begin{proof}
By completing $K$ with respect to the valuation $v_K$ we can assume, without loss of generality, that $K$ is complete.
Let $E/K$ be an elliptic curve with a $K$-rational point of order $5$. Then $E/K$ can be given by a Weierstrass equation of the form 
$$\quad y^2+(1-c)xy-by=x^3-bx^2$$
(see \cite{hus} Section $4.4$), with 
$$b=\lambda, \quad c=\lambda, \quad \text{and} \quad \lambda \in K.$$
Let $\lambda=\frac{s}{t}$ with $s,t \in \mathcal{O}_K$ and such that $v_K(s)=0$ or $v_K(t)=0$.
Using the transformation
$x \rightarrow{x/t^2}$ and  $y \rightarrow{y/t^3}$, we obtain a new Weierstrass equation of the form
$$y^2+ (t-s)xy-st^2y =x^3-stx^2.$$
The discriminant of this Weierstrass equation is $$\Delta=s^5t^5(s^2-11st-t^2),$$
with invariants $$c_4(s,t)=24st^2(-s+t)+(s^2-6st+t^2)^2,$$
and $$c_6(s,t)=-(s^2 + t^2)(s^4 - 18s^3t + 74s^2t^2 + 18st^3 + t^4).$$
Note that if $v_K(s)>0$ or $v_K(t)>0$, then $E/K$ has split multiplicative reduction so in what follows we must have that $v_K(s)=0$ and $v_K(t)=0$ since we assume that $E/K$ has additive reduction.

\begin{claim}\label{claim1wolfram}
There exist polynomials $A(s,t),B(s,t) \in \mathbb{Z}[s,t]$ such that $$A(s,t)c_4(s,t)+B(s,t)c_6(s,t)=2^{12}3^65t^9.$$
\end{claim}

Consider first the polynomials $$f(x)=c_4(x,1)=x^4-12x^3+14x^2+12x+1$$ and $$g(x)=c_6(x,1)=- (x^2+1)(x^4-18x^3+74x^2+18x+1).$$ Using the WOLFRAM \cite{mathematica} command \texttt{PolynomialExtendedGCD[\textit{f}(x),\textit{g}(x)]} we obtain that $f(x)$ and $g(x)$ are coprime in $\mathbb{Q}[x]$ and moreover if $$a(x):=\frac{1}{2^{12}3^65}(6471756x^5-1171065600x^4+4965235200x^3-472780800x^2+4900020480x+698035968)$$ and  $$b(x):=\frac{1}{2^{12}3^65}(64717056x^3-782763264x^2+980543232x+683106048),$$ then $$a(x)f(x)+b(x)g(x)=1.$$
Setting $x=\frac{s}{t}$ and multiplying the last expression by $2^{12}3^65t^9$ we obtain that if $$A(s,t):=a(\frac{s}{t})2^{12}3^65t^5 \; \text{and} \; B(s,t):=b(\frac{s}{t})2^{12}3^65t^3,$$ then $A(s,t)$ and $B(s,t)$ are in $\mathbb{Z}[s,t]$, and $$A(s,t)c_4(s,t)+B(s,t)c_6(s,t)=2^{12}3^65t^9.$$
This proves the claim.

\textit{Proof of $(i)$:} Recall that we assume that $v_K(s)=v_K(t)=0$. If $v_K(c_4(s,t))=0$ or $v_K(c_6(s,t))=0$, then by Tableau 1 of \cite{pap} we find that $E/K$ has semi-stable reduction. If $v_K(c_4(s,t))>0$ and $v_K(c_6(s,t))>0$, then since $v_K(t)=0$ and $v_K(5)=1$, by Claim \ref{claim1wolfram} we obtain that min$\{v_K(c_4), v_K(c_6)\} \leq 1$. As a result, Tableau 1 of \cite{pap} implies that if $E/K$ has additive reduction, then the reduction type can only be II or III.

\textit{Proof of $(ii)$:} Assume that $v_K(5)=2$. Since $v_K(t)=0$, Claim \ref{claim1wolfram} implies that min$\{v_K(c_4), v_K(c_6)\} \leq 2$. Therefore, by Tableau 1 of \cite{pap} we find that if $E/K$ has additive reduction, then the reduction type can only be II, III, IV or I$_n^*$ for some $n \geq 0$.

\textit{Proof of $(iii)$:} Assume now that $v_K(5)=3$. By Claim \ref{claim1wolfram}, since $v_K(t)=0$ and $v_K(5)=3$, we obtain min$\{v_K(c_4), v_K(c_6)\} \leq 3$ so by Tableau 1 of \cite{pap} we find that reduction type II$^*$ cannot occur. 

\textit{Proof of $(iv)$:} Assume that $E/K$ has additive reduction and let $P \in E(K)[10]$. Since $2P$ has order $5$, Part $(i)$ implies that $E/K$ has reduction type II or III. We will show that reduction type II cannot occur. Indeed, consider the point $Q=5P$ which is a point of order $2$. By Theorem $5.9.4$ of \cite{Krummthesis} we get that $Q$ is not in the kernel of the reduction. This implies that $Q$ reduces to the singular point. Indeed, if $Q$ reduces to a non-singular point, then the reduction of $Q$ has order $2$, which is impossible since $\Tilde{E}_{ns}(k) \cong \mathbb{G}_a(k)$ and the latter does not have points of order $2$. Therefore, $Q$ reduces to the singular point. As a result, $2$ divides the component group of $E/\mathbb{Q}$ at $5$. Finally, by the table on page $46$ of \cite{tatealgorithm} we find that reduction type II cannot occur.
\end{proof}

\begin{remark}
We note that if $v_K(5)=1$, then both of the allowed reduction types in Part $(i)$ of Proposition \ref{5torsionprop} occur. For instance, the elliptic curves $E/\mathbb{Q}$ with Cremona labels $50b1$ and $50b2$ have a $\mathbb{Q}$-rational point of order $5$ and reduction of type II modulo $5$. The elliptic curve with Cremona label $150a4$ has a $\mathbb{Q}$-rational point of order $10$ (and hence of order $5$) and reduction of type III modulo $5$.

Concerning Part $(ii)$ of Proposition \ref{5torsionprop}, we first note that if an elliptic curve $E/\mathbb{Q}$ has reduction type II (resp. III) modulo $5$, then after a quadratic extension ramified above $5$ the base change will have reduction type IV (resp. I$_0^*$) modulo the prime above $5$. Moreover, using the LMFDB database we found examples of elliptic curves defined over $\mathbb{Q}(\sqrt{5})$ that have reduction I$_1^*$ modulo the prime above $5$. For instance, the elliptic curve with LMFDB label $2.2.5.1$-$2525.2$-g$1$ is one such example. Thus, the types IV, I$_0^*$, I$_1^*$ in Part $(ii)$ of Proposition \ref{5torsionprop} occur.
\end{remark}

\section{Difference between the mixed and the equicharacteristic cases}\label{section5}

Let $E/K$ be an elliptic curve with a $K$-rational point of order $p$. The reduction types that can occur in the case where $K$ is the fraction field of a discrete valuation ring of characteristic $p \geq 5$ have already been studied by Liedtke and Schr\"{o}er in \cite{liedtkeschroer}. More precisely, among other results in \cite{liedtkeschroer}, they prove the following theorem (see Theorem 4.3 of \cite{liedtkeschroer} and the Theorem in page 2158 of \cite{liedtkeschroer}).

\begin{theorem}
 Let $R$ be a henselian discrete valuation ring of characteristic $p \geq 5$, with fraction field $K$ and residue field $k$ which is assumed to be algebraically closed. Let $E/K$ be elliptic curve with additive reduction and containing a $K$-rational point of order $p$. Then $E/K$ has potentially supersingular reduction. Moreover,
 \begin{enumerate}[topsep=2pt,label=(\roman*)]
\itemsep0em 
\item If $p \equiv 1 \; (mod \; 12)$, then $E/K$ can only have reduction type I$_0^*$.
\item If $p \equiv 5 \; (mod \; 12)$, then $E/K$ can only have reduction type II, IV, I$_0^*$, IV$^*$ or II$^*$.
\item If $p \equiv 7 \; (mod \; 12)$, then $E/K$ can only have reduction type III, I$_0^*$ or III$^*$.
\item If $p \equiv 11 \; (mod \; 12)$, then $E/K$ can only have reduction type II, III, IV, I$_0^*$, IV$^*$, III$^*$ or II$^*$.
\end{enumerate}
\end{theorem}

We note that it is not true in the mixed characteristic case that an elliptic curve $E/K$ with additive reduction and integral $j$-invariant has potentially {\it supersingular} reduction when it has a $K$-rational point of prime order $p \geq 5$. There are examples of elliptic curves defined over $\mathbb{Q}$ with a $7$-torsion point that have additive reduction of Kodaira type II (see Remark \ref{remark7torsion}). The following proposition shows that for $p=7$ the reduction of such elliptic curves cannot be potentially supersingular.

\begin{proposition}\label{proppotsupersingular}
  Let $R$ be a discrete valuation ring with fraction field $K$, where $K$ is of any characteristic, and with residue field characteristic $p \geq 5$. Let $E$ be an elliptic curve over $K$. Then we have the following:
   \begin{enumerate}[topsep=2pt,label=(\roman*)]
\itemsep0em 
    \item  If $E/K$ has reduction of type II, II$^*$, IV or IV$^*$, then the reduction is potentially supersingular if and only if $p \equiv 2 \; (\text{mod}\; 3)$.
      \item If $E/K$ has reduction of type III or III$^*$,  then the reduction is potentially supersingular if and only if $p \equiv 3 \; (\text{mod}\; 4)$.
\end{enumerate}
 \end{proposition}

 \begin{proof}
  Assume that $E/K$ has reduction type II, II$^*$, IV or IV$^*$. By the table on page $46$ of \cite{tatealgorithm} the $j$-invariant of $E/K$ reduces to zero in the residue field. However, $0$ is a supersingular $j$-invariant if and only if  $p \equiv 2 \; (\text{mod}\; 3)$ by \cite{aec} Example V.4.4. This proves Part $(i)$.
 
 Assume that $E/K$ has reduction type III or III$^*$. By the table on page $46$ of \cite{tatealgorithm} the $j$-invariant of $E/K$ reduces to $1728$ in the residue field. Since $1728$ is a supersingular $j$-invariant if and only if  $p \equiv 3 \; (\text{mod}\; 4)$ by \cite{aec} Example V.4.4, we have that reduction is potentially supersingular if and only $p \equiv 3 \; (\text{mod}\; 4)$. This proves Part $(ii)$.
 \end{proof}

\bibliographystyle{plain}
\bibliography{bibliography.bib}

\end{document}